\documentclass[10pt]{article}
\usepackage[T1]{fontenc}
\usepackage[utf8]{inputenc}

\usepackage{lmodern}
\usepackage{xspace}
\usepackage[normalem]{ulem}
\usepackage{amsfonts,amsmath,amssymb, amsthm, mathtools}
\usepackage{thm-restate}
\usepackage{dsfont} 
\usepackage{algorithmicx,algpseudocode,algorithm}

\usepackage[usenames,dvipsnames,table]{xcolor}

\usepackage{csquotes}

\usepackage{relsize}

\usepackage{multirow}
\usepackage{chngpage} 

\usepackage{mfirstuc}

\usepackage{tikz}
\usetikzlibrary{arrows}
\usetikzlibrary{calc,decorations.pathmorphing,patterns}
\usepackage[colorlinks,citecolor=blue,bookmarks=true,linktocpage]{hyperref}
\usepackage{aliascnt}
\usepackage[numbered]{bookmark}

\usepackage{fullpage}

\usepackage{titling}

\usepackage[shortlabels]{enumitem}
  \setitemize{noitemsep,topsep=3pt,parsep=2pt,partopsep=2pt} \setenumerate{itemsep=1pt,topsep=2pt,parsep=2pt,partopsep=2pt}
  \setdescription{itemsep=1pt}

\usepackage{verbatim}

\usepackage{mleftright}  \makeatletter
\@ifundefined{theorem}{\theoremstyle{plain} \newtheorem{theorem}{Theorem}\newaliascnt{coro}{theorem}
  	  \newtheorem{corollary}[coro]{Corollary}
  	\aliascntresetthe{coro}
  	\newaliascnt{lem}{theorem}
  		\newtheorem{lemma}[lem]{Lemma}
  	\aliascntresetthe{lem}
  	\newaliascnt{clm}{theorem}
  		
	\aliascntresetthe{clm}
	\newaliascnt{factcnt}{theorem}
 	 	\newtheorem{fact}[theorem]{Fact}
	\aliascntresetthe{factcnt}
  	
  \newaliascnt{prop}{theorem}
  		
	\aliascntresetthe{prop}
	\newaliascnt{conj}{theorem}
  		
	\aliascntresetthe{conj}
  \theoremstyle{remark} 
  \newtheorem{remark}[theorem]{Remark}
  	\newtheorem{question}[theorem]{Question}

  \theoremstyle{definition} \newaliascnt{defn}{theorem}
 		 
 	 \aliascntresetthe{defn}
}{}
\makeatother

\providecommand{\email}[1]{\href{mailto:#1}{\nolinkurl{#1}\xspace}}
\newcommand{\ignore}[1]{\leavevmode\unskip}  
\newcommand{\eps}{\ensuremath{\varepsilon}\xspace}
     \newcommand{\eqdef}{\stackrel{\rm def}{=}}

\newcommand{\domain}{\ensuremath{\Omega}\xspace}

\newcommand{\bigOmega}[1]{{\Omega\mleft( #1 \mright)}}

 			\newcommand{\indicSet}[1]{\mathds{1}_{#1}}                                                                                           
\newcommand{\dtv}{\operatorname{d}_{\rm TV}}
\newcommand{\kl}{\operatorname{KL}}

\newcommand{\kldiv}[2]{{\kl\mleft({#1 \,\|\, #2}\mright)}}

\newcommand{\totalvardistrestr}[3][]{{\dtv^{#1}\mleft({#2, #3}\mright)}}
\newcommand{\totalvardist}[2]{\totalvardistrestr[]{#1}{#2}}

\newcommand\restr[2]{{\left.\kern-\nulldelimiterspace #1 \vphantom{\big|} \right|_{#2} }}

\newcommand{\shortexpect}{\mathbb{E}}

\newcommand{\bernoulli}[1]{\ensuremath{\operatorname{Bern}\!\left( #1 \right)}}

\newcommand{\norm}[1]{\lVert#1{\rVert}}
\newcommand{\normone}[1]{{\norm{#1}}_1}

\newcommand{\norminf}[1]{{\norm{#1}}_\infty}
\newcommand{\abs}[1]{\left\lvert #1 \right\rvert}

\newcommand{\R}{\ensuremath{\mathbb{R}}\xspace}

\newcommand{\pdfsamp}{dual\xspace}
\newcommand{\cdfsamp}{cumulative dual\xspace}
\newcommand{\Pdfsamp}{\expandafter\capitalisewords\expandafter{\pdfsamp}}
\newcommand{\Cdfsamp}{\expandafter\capitalisewords\expandafter{\cdfsamp}}

\makeatletter

\newcommand{\Rom}[1]{\expandafter\@slowromancap\romannumeral #1@}

\makeatother

\makeatletter
  \AtBeginDocument{
  \begingroup
  \toks@={}\toksdef\toks@@=2 \toks@@={}\long\def\@ReturnFiFi#1#2\fi\fi{\fi\fi#1}\def\scan@author#1#2 \and#3\@nil{\ifx\\#3\\\ifcase#1 \toks@={#2}\else
      \ifnum#1>1 \toks@=\expandafter{\the\expandafter\toks@\expandafter,\expandafter\space
          \the\toks@@
        }\fi
      \toks@=\expandafter{\the\toks@\space and #2}\fi
    \else
      \ifcase#1 \toks@={#2}\@ReturnFiFi{\scan@author1#3\@nil
        }\else
        \ifnum#1>1 \toks@=\expandafter{\the\expandafter\toks@\expandafter,\expandafter\space
            \the\toks@@
          }\fi
      \toks@@={#2}\@ReturnFiFi{\scan@author2#3\@nil
      }\fi
  \fi
  }\expandafter\expandafter\expandafter\scan@author
  \expandafter\expandafter\expandafter0\expandafter\@author\space\and\@nil
  \edef\x{\endgroup
  \noexpand\hypersetup{pdfauthor={\the\toks@}}}\x
  }
\makeatother
 \usepackage{pgfplots}
\usepgfplotslibrary{fillbetween}
\pgfplotsset{compat=1.17} 
\usepackage{framed}

\newcommand{\ns}{n}

\renewcommand{\eqdef}{\coloneqq}
\usepackage{utopia}
\usepackage{filecontents}

\title{A short note on an inequality between KL and TV}
\author{Cl\'ement L. Canonne\thanks{University of Sydney. Email: \email{clement.canonne@sydney.edu.au}}}
\newcommand{\p}{\mathbf{p}}
\newcommand{\q}{\mathbf{q}}

\begin{document}
\maketitle

The goal of this short note is to discuss the relation between Kullback--Leibler divergence and total variation distance, starting with the celebrated Pinsker's inequality relating the two, before switching to a simple, yet (arguably) more useful inequality, apparently not as well known. We summarize the gist of it below:

\begin{theorem}[The BH Bound]
  \label{theo:bh}
    For every two probability distributions $\p,\q$, we have the simple yet never vacuous bound
    \begin{equation}
        \totalvardist{\p}{\q} \leq \sqrt{1-e^{-\kldiv{\p}{\q}}}
    \end{equation}
    and this is never worse than Pinsker's inequality (except for a factor $\sqrt{2}$ for $\kldiv{\p}{\q}\ll 1$.)
\end{theorem}
While establishing this bound and discussing its various aspects, we will consider probability distributions over some set $\domain$, and conveniently ignore any measurability or absolute continuity issue~--~the reader is encouraged to think of discrete $\domain$ for concreteness. Everything does however apply to the general setting, given the suitable insertion of the words ``Radon--Nikodym derivative'' and ``measurable'' in appropriate locations.

\paragraph{Total variation distance and Kullback--Leibler divergence.}

The TV distance and KL divergence (in nats) between two probability distributions $\p,\q$ over $\domain$ are given respectively by
\[
    \totalvardist{\p}{\q} = \sup_{S\subseteq \domain} (\p(S)-\q(S)) = \frac{1}{2}\normone{\p-\q}\in[0,1]
\]
and
\[
    \kldiv{\p}{\q}  = \sum_{x\in\domain} \p(x)\log\frac{\p(x)}{\q(x)} \in [0,\infty)
\]
where $\log$ is the natural logarithm, with the convention that $0\log 0 = 0$. Both TV distance and KL divergence are special cases of what is known as \emph{$f$-divergences}, and they both enjoy a lot of crucial properties, such as the data processing inequality, which we will not get into here.\footnote{The KL divergence, annoyingly, is not symmetric in its arguments and thus not a real metric, but it makes up for it sometimes.}

\paragraph{Organisation.}
We start by a (very brief) review of Pinsker's inequality, and its shortcomings, in~\autoref{sec:pinsker}, before stating and deriving the BH bound in~\autoref{sec:bh}. The reader asking themselves why we should care at all about this improved bound can skip directly to~\autoref{sec:applications} for some motivation and applications, and those keen on the Donsker--Varadhan formula (or looking for an open question) might enjoy~\autoref{sec:tfl}. Finally,~\autoref{sec:discussion} provides some pointers, and discusses a slightly more refined (albeit much more unwieldy) bounds.

\section{Pinsker's inequality}
  \label{sec:pinsker}
  
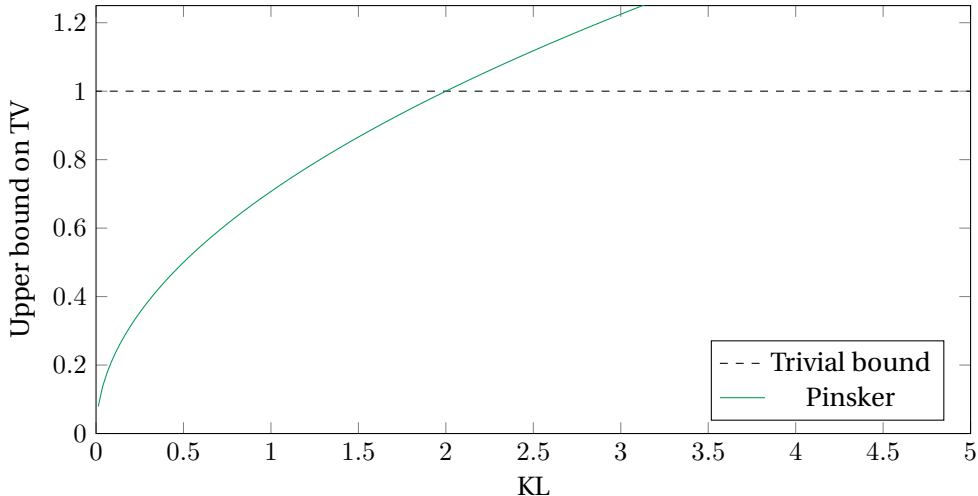
\begin{figure}[H]\centering
\begin{tikzpicture}
  \begin{axis}[ 
    xlabel={KL},
    ylabel={Upper bound on TV},
    xmin=0, xmax=5, ymin=0, ymax=1.25,
    width=0.8\textwidth,
    height=\axisdefaultheight,
    legend pos=south east,
  ] 
    \addlegendentry{Trivial bound}
    \addlegendentry{Pinsker}
    \addplot[samples=100, color=black, dashed] {1}; 
    \addplot[samples=400, color=ForestGreen] {sqrt(x/2)};
  \end{axis}
\end{tikzpicture} 
\caption{\label{fig:pinsker}Pinsker's inequality becomes vacuous for $\kldiv{\p}{\q}>2$. That's a downer.}
\end{figure}

  We first state our baseline, Pinsker's inequality, a fundamental relation between KL divergence and total variation distance originally due to, well, Pinsker~\cite{Pinsker64}, although in a weaker form and with suboptimal constants: the constant was then independently improved to the optimal $1/\sqrt{2}$ by Kullback, Csizs\'ar, and Kemperman. See~\cite[Section~2.8]{Tsybakov09} for a discussion.
\begin{lemma}[Pinsker's Inequality]
  \label{lemma:pinsker}
For every $\p,\q$ on $\domain$, 
\begin{equation}
  \label{eq:pinsker}
  \totalvardist{\p}{\q} \leq \sqrt{\frac{1}{2}\kldiv{\p}{\q}}\,.
\end{equation}
\end{lemma}
There are many proofs of Pinsker's inequality: e.g., \cite[Lemma~2.5]{Tsybakov09}, or a very clever argument due to Pollard, or even in~\autoref{sec:tfl} of this very note, using the Donsker--Varadhan formula (also known by me as ``Thomas' Favourite Lemma''). One can for instance consult~\cite{IlyaRaz:MO} for a list; we will here follow an argument from Yihong Wu's lecture notes~\cite[Theorem~4.5]{Wu20}, which has the advantage of being nearly magical.
\begin{proof}
Consider first the binary case, i.e., where $\p$,$\q$ are Bernoulli distributions $\bernoulli{p}$ and $\bernoulli{q}$, respectively. Then 
$
    \totalvardist{\p}{\q} = \abs{p-q}
$
and so we are left to prove
\begin{equation}
  \label{eq:proving:pinsker:binary}
    2(p-q)^2 \leq p\log\frac{p}{q} + (1-p)\log\frac{1-p}{1-q}
\end{equation}
Note that the cases where either $p$ or $q$ is in $\{0,1\}$ are easily checked \textit{(verify it!)}, so we can assume $p,q\in(0,1)$. 
To prove~\eqref{eq:proving:pinsker:binary} in this case, we introduce the function $f\colon(0,1)\to\R$ defined by $f(x) = p\log x + (1-p)\log(1-x)$, and observe that the RHS of~\eqref{eq:proving:pinsker:binary} is exactly $f(p)-f(q)$. We then can write
\[
    f(p)-f(q) = \int_{q}^p f'(x)dx = \int_{q}^p \frac{p-x}{x(1-x)}dx \geq 4\int_{q}^p (p-x)dx = 4\cdot \frac{1}{2}(p-q)^2
\]
establishing~\eqref{eq:proving:pinsker:binary} (note that we used the fact that $x(1-x) \leq 1/4$ for $x\in(0,1)$).\smallskip

Turning to the general case, let $\p,\q$ be distributions on an arbitrary domain $\domain$, and fix any measurable subset $S\subseteq \domain$. For $X,Y$ distributed according to $\p$ and $\q$, the random variables $\indicSet{S}(X)$ and $\indicSet{S}(Y)$ are have distributions $\p'\eqdef \bernoulli{\p(S)}$ and $\q' \eqdef \bernoulli{\q(S)}$ respectively, and therefore
\[
    2(\p(S)-\q(S))^2 = 2\totalvardist{\p'}{\q'}^2 \leq \kldiv{\p'}{\q'} \leq \kldiv{\p}{\q}
\]
where the first inequality is~\eqref{eq:proving:pinsker:binary}, and the second is the data processing inequality. Since this inequality holds for every $S$, taking a supremum over $S$ leads to
\[
    2\totalvardist{\p}{\q}^2 = 2\sup_{S\subseteq\domain}(\p(S)-\q(S))^2 \leq \kldiv{\p}{\q}\,,
\]
establishing Pinsker's inequality.
\end{proof}
Before we try to improve upon Pinsker's inequality, let us note that one particular avenue is doomed: specifically, the constant $1/\sqrt{2}$ in~\eqref{eq:pinsker} cannot be replaced by any $c<1/\sqrt{2}$. To see why, fix any $\eps\in(0,1/4)$, and observe that for $\p=\bernoulli{1/2}$ and $\q=\bernoulli{1/2+\eps}$ we have $\totalvardist{\p}{\q} = \eps$ and $\kldiv{\p}{\q}=\frac{1}{2}\log\frac{1}{1-4\eps^2}$, so that
\begin{equation}
  \frac{\kldiv{\p}{\q}}{\totalvardist{\p}{\q}^2} = -\frac{\log(1-4\eps^2)}{2\eps^2} \xrightarrow[\eps\to 0^+]{}2\,.
\end{equation}
Still, in spite of its multiple applications in Statistics and information theory and its ``optimality'' shown above, Pinsker's inequality suffers a major drawback: by definition, the TV distance is always at most $1$, yet the RHS of~\eqref{eq:pinsker} grows unbounded with the KL divergence. In other terms, the bound is totally and utterly useless for any $\kldiv{\p}{\q}>2$, as depicted in~\autoref{fig:pinsker}.
To see why one would care about this issue (without jumping yet to~\autoref{sec:applications}), consider the following very simple and intuitive fact: \emph{``if $\kldiv{\p}{\q}<\infty$, then $\totalvardist{\p}{\q} < 1$.''}  While absolutely true, this claim cannot be proven from Pinsker's inequality. Even worse, using Pinsker's one cannot even establish that if $\kldiv{\p}{\q}< 2.01$ then the two distributions $\p$ and $\q$ have TV distance bounded away from $1$!

\section{The Bretagnolle--Huber bound}
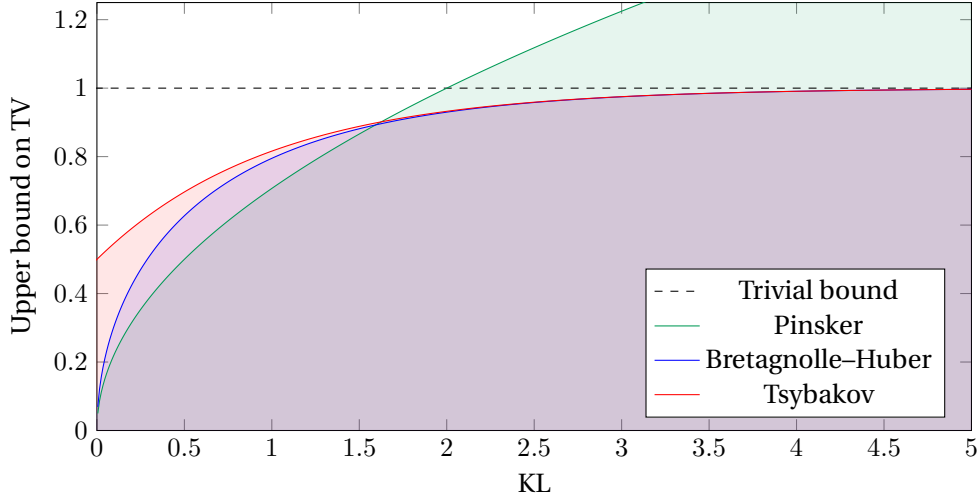
\begin{figure}[H]\centering
\begin{tikzpicture}
  \begin{axis}[ 
    xlabel={KL},
    ylabel={Upper bound on TV},
    xmin=0, xmax=5, ymin=0, ymax=1.25,
    width=0.8\textwidth,
    height=\axisdefaultheight,
    legend pos=south east,
  ] 
    \addlegendentry{Trivial bound}
    \addlegendentry{Pinsker}
    \addlegendentry{Bretagnolle--Huber}
    \addlegendentry{Tsybakov}
    \addplot[samples=5, color=black, dashed] {1}; 
    \addplot[samples=1000, color=ForestGreen, name path=Pinsker] {sqrt(x/2)};
    \addplot[samples=1000, color=blue, name path=BH] {sqrt(1-e^(-x))}; 
    \addplot[samples=1000, color=red, name path=Tsybakov] {1-1/2*e^(-x)}; 
    
    \addplot[draw=none,mark=none,name path=Baseline] {0};     \addplot[ForestGreen,fill opacity=0.1] fill between[of=Pinsker and Baseline];
    \addplot[blue,fill opacity=0.1] fill between[of=BH and Baseline];
    \addplot[red,fill opacity=0.1] fill between[of=Tsybakov and Baseline];
  \end{axis}
\end{tikzpicture} 
\caption{\label{fig:bounds:pinsker:bh:tsybakov}The four upper bounds we have: shaded regions correspond to values of TV still allowed by the corresponding bounds (so smaller shaded areas are better). As we can see, Pinsker's bound~\eqref{eq:pinsker} is useful for small values of $\kldiv{\p}{\q}$ only. Tsybakov's bound~\eqref{eq:bh:tsybakov} is much better for large $\kldiv{\p}{\q}$, and in particular have the right behaviour as $\kldiv{\p}{\q}\to\infty$, never becoming worse than the trivial bound. However, it is now useless for small $\kldiv{\p}{\q}$, and does not even go to $0$ as $\kldiv{\p}{\q}\to 0^+$. The clear winner is the Bretagnolle--Huber bound~\eqref{eq:bh}, which not only is never worse than Tsybakov's (obviously), but also has the right behaviour for small values of $\kldiv{\p}{\q}$, being essentially equivalent (up to a constant factor) to Pinsker's in that regime.}
\end{figure}

In view of the above, can we hope for some better inequality which does not leap into vacuousness when the KL divergence gets large? The answer is, thankfully, yes.
  \label{sec:bh}
\begin{lemma}[The BH Bound]
  \label{lemma:bh}
For every $\p,\q$ on $\domain$, 
\begin{equation}
  \label{eq:bh}
  \totalvardist{\p}{\q} \leq \sqrt{1-e^{-\kldiv{\p}{\q}}}\,.
\end{equation}
\end{lemma}
\begin{proof}
We follow the original argument of~\cite[Lemma~2.1]{BretagnolleH78}: fixing $\p,\q$, we define, for $X$ distributed according to $\p$, the random variables $U\eqdef \frac{\q(X)}{\p(X)}$, $V\eqdef (U-1)_+$, and $W\eqdef 1+V-U=(1-U)_+$. One can check that
\[
    \totalvardist{\p}{\q} = \frac{1}{2}\shortexpect_{\p}[\abs{U-1}] = \shortexpect_{\p}[V] = \shortexpect_{\p}[W]
\]
and that by construction $(1+V)(1-W)=U$, so that $\log U = \log(1+V)+\log(1-W)$. Moreover, since
$
  \kldiv{\p}{\q} = -\sum_{x\in\domain}\p(x)\log\frac{\q(x)}{\p(x)} = -\shortexpect_{\p}[\log U]
$, we get by Jensen's inequality that
\[
  -\kldiv{\p}{\q} = \shortexpect_{\p}[\log(1+V)]+\shortexpect_{\p}[\log(1-W)] \leq \log(1+\shortexpect_{\p}[V])+\log(1-\shortexpect_{\p}[W])
= \log(1-\totalvardist{\p}{\q}^2)
\]
which, exponentiating both sides, rearranging and taking the square root, proves the lemma.
\end{proof}
Now, instead of the above bound, one may encounter the following weaker one, for instance in Tsybakov's monograph~\cite{Tsybakov09}.\footnote{which is worth reading by itself, as it countains many gems, insights, and useful discussions.} It is unclear to me what advantage this looser inequality holds over~\eqref{lemma:bh}, but as we shall see in~\autoref{fig:bounds:pinsker:bh:tsybakov} it at least behaves in a satisfying way for large values of $\kldiv{\p}{\q}$, and is never bigger than the trivial bound of $1$.
\begin{corollary}[Tsybakov's version]
  \label{lemma:bh:tsybakov}
For every $\p,\q$ on $\domain$, 
\begin{equation}
  \label{eq:bh:tsybakov}
  \totalvardist{\p}{\q} \leq 1-\frac{1}{2}e^{-\kldiv{\p}{\q}}\,.
\end{equation}
\end{corollary}
\begin{proof}
  This readily frollows from~\autoref{lemma:bh}, upon noting that $\sqrt{1-e^{-x}} \leq 1-\frac{1}{2}e^{-x}$ for all $x\in[0,\infty)$ (just square both sides and expand the RHS).
\end{proof}

Because it is somewhat fun to do, we also reproduce Tsybakov's proof of~\autoref{lemma:bh:tsybakov}, and show how it can be used to derive~\autoref{lemma:bh} (so I am at a loss as to why Tsybakov would only state the weaker version in his monograph).\footnote{We observe that the same ``extension'' of Tsybakov's argument can be found in the proof of~\cite[Lemma 6]{GaoHRZ19}.}
\begin{proof}[Proof of~\autoref{lemma:bh} from the argument of~\cite{Tsybakov09}, Lemma~2.6]
Fix $\p,\q$. First, we observe that one can write 
\begin{equation}
  \totalvardist{\p}{\q}=1-\sum_{x\in\domain} \min(\p(x),\q(x)) = \sum_{x\in\domain} \max(\p(x),\q(x)) - 1
\end{equation}
\emph{(this is a useful trick, check it!)}, and therefore by Cauchy--Schwarz
\begin{align}
    1-\totalvardist{\p}{\q}^2 
    &= (1+\totalvardist{\p}{\q})(1-\totalvardist{\p}{\q})
    = \mleft(\sum_{x\in\domain} \max(\p(x),\q(x))\mright)\mleft(\sum_{x\in\domain} \min(\p(x),\q(x))\mright) \notag\\
    &\geq \mleft(\sum_{x\in\domain} \sqrt{\max(\p(x),\q(x))\min(\p(x),\q(x))}\mright)^2
    = \mleft(\sum_{x\in\domain} \sqrt{\p(x)\q(x)}\mright)^2 \label{eq:tsybakov:step1}
\end{align}
which will come handy very soon. Indeed, what Tsybakov does show is the following:\footnote{The eagle-eyed reader may recognize in the LHS the square of the \emph{Hellinger affinity}: indeed, this inequality actually provides a bound on the Hellinger distance in terms of KL divergence, which in turn implies the bound on TV distance by~\eqref{eq:tsybakov:step1}.}
\begin{align}
    \mleft(\sum_{x\in\domain} \sqrt{\p(x)\q(x)}\mright)^2
    &= e^{2\log \sum_{x} \sqrt{\p(x)\q(x)}}
    = e^{2\log \sum_{x} \p(x)\sqrt{\frac{\q(x)}{\p(x)}}} \notag\\
    &= e^{2\log \shortexpect_\p\mleft[\sqrt{\frac{\q(X)}{\p(X)}}\mright]} 
    \geq e^{2\shortexpect_\p\mleft[\log \sqrt{\frac{\q(X)}{\p(X)}}\mright]} \tag{Jensen's inequality}\notag\\
    &= e^{\shortexpect_\p\mleft[\log \frac{\q(X)}{\p(X)}\mright]} 
    = e^{-\kldiv{\p}{\q}} \label{eq:tsybakov:step2}
\end{align}
(to be precise, sums and expectations are restricted to the support of $\p$, to avoid dividing by zero). Combining~\eqref{eq:tsybakov:step1} and~\eqref{eq:tsybakov:step2} yields~\autoref{lemma:bh}.
\end{proof}

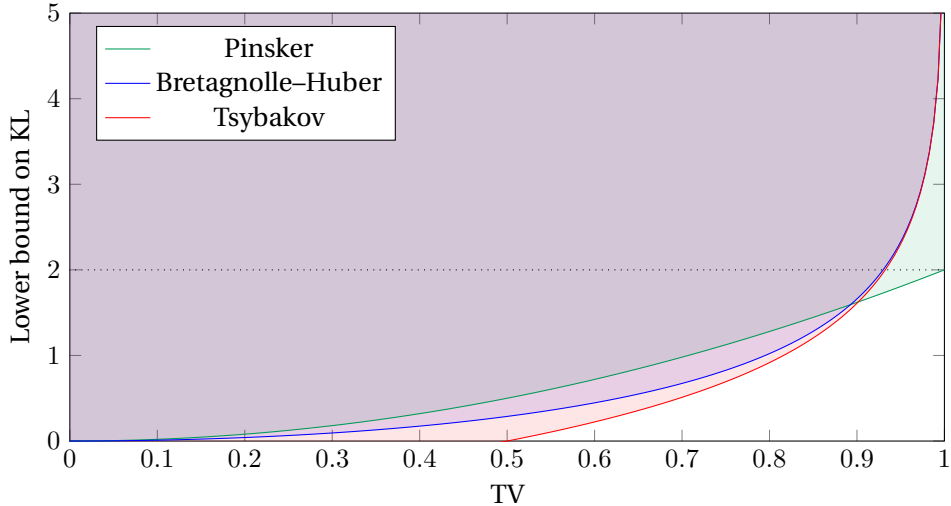
\begin{figure}[ht!]\centering
\begin{tikzpicture}
  \begin{axis}[ 
    xlabel={TV},
    ylabel={Lower bound on KL},
    xmin=0, xmax=1, ymin=0, ymax=5,
    width=0.8\textwidth,
    height=\axisdefaultheight,
    legend pos=north west,
  ] 
\addlegendentry{Pinsker}
    \addlegendentry{Bretagnolle--Huber}
    \addlegendentry{Tsybakov}
    \addplot[samples=500, color=ForestGreen,name path=Pinsker] {2*x*x};
    \addplot[samples=2000, color=blue,name path=BH] {-ln(1-x*x)}; 
    \addplot[samples=2000, color=red,name path=Tsybakov] {-ln(2*(1-x))}; 
    \addplot[samples=5, color=black, dotted] {2}; 
    
    \addplot[draw=none,mark=none,name path=Baseline] {5};     \addplot[ForestGreen,fill opacity=0.1] fill between[of=Pinsker and Baseline];
    \addplot[blue,fill opacity=0.1] fill between[of=BH and Baseline];
    \addplot[red,fill opacity=0.1] fill between[of=Tsybakov and Baseline];
  \end{axis}
\end{tikzpicture} 
\caption{\label{fig:bounds:pinsker:bh:tsybakov:reciprocal}The lower bounds which~\eqref{eq:pinsker},~\eqref{eq:bh}, and~\eqref{eq:bh:tsybakov} give on the KL divergence. The shaded areas are the values of KL (as a function of TV) still allowed by the corresponding inqualities, so smaller shaded area is better: as one can see, Pinsker's inequality is unable to rule out any value of KL greater than 2, while the bound given by Tsybakov only kicks in for $\text{TV}\geq 1/2$.}
\end{figure}

\noindent To see how these bounds \eqref{eq:pinsker}, \eqref{eq:bh}, and \eqref{eq:bh:tsybakov} compare, let us look at a plot.
\autoref{fig:bounds:pinsker:bh:tsybakov} clearly hints that the BH bound obtained in~\autoref{lemma:bh} is never much worse than the one from Pinsker's inequality, but let us make this formal. First, a Taylor approximation shows that $\sqrt{1-e^{-x}} = \sqrt{x} + o(\sqrt{x})$ as $x\to 0^+$, so for small TV our new bound is worse than Pinsker's by only a factor $\sqrt{2}$. It is actually easy to see that this is always the case, as the inequality $\sqrt{1-e^{-x}} \leq \sqrt{2}\cdot \sqrt{\frac{x}{2}}$ (for $x\geq 0$) is equivalent  to $1-x\leq e^{-x}$, which holds by convexity. We can summarize this as follows:
\begin{framed}
\noindent The BH bound~\eqref{eq:bh} is never vacuous, has the right behaviour when $\kldiv{\p}{\q}\to\infty$ \emph{and} $\kldiv{\p}{\q}\to 0^+$, and is at worst a $\sqrt{2}$ factor off from Pinsker's bound~\eqref{eq:pinsker}.
\end{framed}
To provide a complementary view of the depiction of the bounds from~\autoref{fig:bounds:pinsker:bh:tsybakov} (which showed the \emph{upper} bounds on TV, as a function of KL, implied by the Pinsker, BH, and Tsybakov inequalities), we give in~\autoref{fig:bounds:pinsker:bh:tsybakov:reciprocal} the corresponding \emph{lower} bounds on KL as a function of TV.

\begin{figure}[ht!]\centering
\begin{tikzpicture}
  \begin{axis}[ 
    xlabel={KL},
    ylabel={Upper bound on TV},
    xmin=0, xmax=5, ymin=0, ymax=1.25,
    width=0.8\textwidth,
    height=\axisdefaultheight,
    legend pos=south east,
  ] 
    \addlegendentry{Trivial bound}
    \addlegendentry{Pinsker}
    \addlegendentry{Bretagnolle--Huber}
    \addlegendentry{``Weak Bretagnolle--Huber''}
    \addplot[samples=100, color=black, dashed] {1}; 
    \addplot[samples=400, color=ForestGreen] {sqrt(x/2)}; 
    \addplot[samples=400, color=blue] {sqrt(1-e^(-x))}; 
    \addplot[samples=400, color=orange] {sqrt(1/(1-e^(-2)))*sqrt(1-e^(-x))}; 
  \end{axis}
\end{tikzpicture} 
\caption{\label{fig:bh:weak}A ``weak Bretagnolle--Huber bound''~\eqref{eq:bh:weak} can be derived from Pinsker's inequality, but it is not a good idea.}
\end{figure}
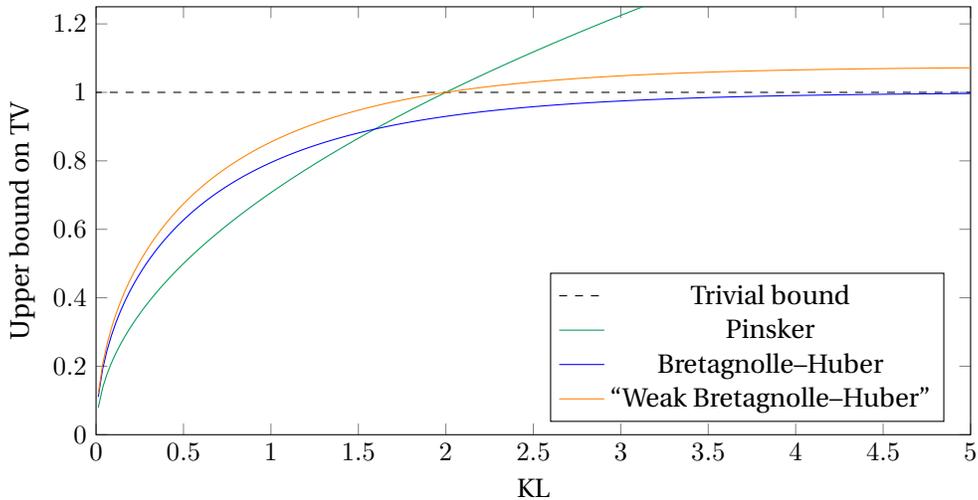

\paragraph{So close, yet so far?}
Interestingly, one can derive an inequality looking similar to the BH bound from Pinsker's inequality, with a major caveat. Note that~\eqref{eq:pinsker} can be equivalently rephrased as follows:
\begin{equation}
  1-e^{-2\totalvardist{\p}{\q}^2} \leq 1-e^{-\kldiv{\p}{\q}}\,.
\end{equation}
Using the (tight) inequality $x\leq (1-e^{-2})^{-1}(1-e^{-2x})$, which holds for all $x\in[0,1]$, we then get
\begin{equation}
  \label{eq:bh:weak}
  \totalvardist{\p}{\q} \leq \frac{1}{\sqrt{1-e^{-2}}}\cdot\sqrt{1-e^{-\kldiv{\p}{\q}}}\,,
\end{equation}
which, except for this leading factor $\frac{1}{\sqrt{1-e^{-2}}} \approx 1.075$, looks very much like~\eqref{eq:bh}. Unfortunately, this leading factor is exactly what makes~\eqref{eq:bh:weak} useless, as the bound is still vacuous whenever $\kldiv{\p}{\q}>2$, and further is strictly weaker than Pinsker's for $\kldiv{\p}{\q}<2$. See~\autoref{fig:bh:weak} for an illustration.

\section{Why do we care?}
  \label{sec:applications}
As we saw, the BH bound~\eqref{eq:bh} and the weaker bound~\eqref{eq:bh:tsybakov} both improve on Pinsker's inequality~\eqref{eq:pinsker} in the regime of large KL. Since for vanishingly small KL Pinsker's inequality is tight in general, and provides a good bound for small constant values as well, it is natural to wonder why we may care about the regime $\operatorname{KL}\gg 1$. One example lies in proving sample complexity lower bounds. As an example and motivation, consider the following (true) fact:
\begin{fact}
  \label{fact:test:bernoulli}
  The number of independent tosses required to distinguish with probability at least $1-\delta$ between a fair coin (i.e., $\bernoulli{1/2}$) and an $\eps$-biased coin (i.e., $\bernoulli{1/2+\eps}$) is $\bigOmega{\log(1/\delta)/\eps^2}$.
\end{fact}
To prove this for constant $\delta>0$, say $\delta=1/10$, the standard way to proceed is to observe that, by a relatively standard argument, we need the number $\ns$ of samples (tosses) to satisfy
\begin{equation}
    \totalvardist{\bernoulli{1/2}^{\otimes\ns}}{\bernoulli{1/2+\eps}^{\otimes\ns}} \geq 1-2\delta
\end{equation}
and we can then use Pinsker's inequality and additivity of KL divergence for product distributions to get
\begin{align*}
    (1-2\delta)^2 
    &\leq \totalvardist{\bernoulli{1/2}^{\otimes\ns}}{\bernoulli{1/2+\eps}^{\otimes\ns}}^2    \\
    &\leq \frac{1}{2}\kldiv{\bernoulli{1/2}^{\otimes\ns}}{\bernoulli{1/2+\eps}^{\otimes\ns}}  \tag{Pinsker}\\
    &= \ns\cdot \frac{1}{2}\kldiv{\bernoulli{1/2}}{\bernoulli{1/2+\eps}}  \\
    &= \ns\cdot \frac{1}{2}\log\frac{1}{1-4\eps^2}  \tag{Direct computation of KL}
\end{align*}
which is at most $4\ns\eps^2$ for $\eps$ small enough, e.g., $0<\eps<1/3$. This shows the $\bigOmega{1/\eps^2}$ for constant $\delta\in(0,1/2)$. It is not hard to see, unfortunately, that this approach will never yield any bound better than $\bigOmega{1/\eps^2}$, even as $\delta\to 0^+$; exactly because Pinsker's inequality does not allow us to discrimate between ``moderately large KL'' and ''KL going to $\infty$.'' But if we were to use the BH bound instead, then the exact same argument shows that we need
\begin{align*}
    (1-2\delta)^2 
    &\leq \totalvardist{\bernoulli{1/2}^{\otimes\ns}}{\bernoulli{1/2+\eps}^{\otimes\ns}}^2    \\
    &\leq 1- e^{-\kldiv{\bernoulli{1/2}^{\otimes\ns}}{\bernoulli{1/2+\eps}^{\otimes\ns}}}  \tag{BH}\\
    &= 1-e^{-\ns\cdot \frac{1}{2}\kldiv{\bernoulli{1/2}}{\bernoulli{1/2+\eps}}}  \\
    &= 1-e^{-\ns\cdot \frac{1}{2}\log\frac{1}{1-4\eps^2}}
\end{align*}
or, reorganizing, $\ns\geq \frac{2}{\log\frac{1}{1-4\eps^2}}\log\frac{1}{1-(1-2\delta)^2} \geq \frac{1}{2\eps^2}\log\frac{1}{2\delta}$ (the last inequality again for $\eps\in(0,1/3)$), which proves~\autoref{fact:test:bernoulli}.

\begin{remark}
  One can also use~\eqref{eq:bh:tsybakov} to prove~\autoref{fact:test:bernoulli} in a similar fashion, which turns out to be even (marginally) simpler: verify it!
\end{remark}

\section{The TFL}
  \label{sec:tfl}
Let us switch gears a little, and consider the relation between these inequalities and the fundamental lemma below, sometimes known as the \emph{Gibbs variational principle}, or the \emph{Donsker--Varadhan formula}~\cite{DonskerV75}, and which I have been told is a special case of Fenchel duality. We will refer to it, succinctly, as \emph{Thomas' Favourite Lemma}, thus named after \href{http://www.thomas-steinke.net/}{Thomas Steinke} and his fondness for this result.
\begin{lemma}[Thomas' Favourite Lemma (TFL)]
  \label{lemma:tfl}
For every $\q\ll\p$, 
\[
    \kldiv{\p}{\q} = \sup_f \mleft( \shortexpect_{\p}[f(X)] - \log \shortexpect_{\q}[e^{f(Y)}] \mright)
\]
where the supremum is over all (measurable) $f\colon\domain\to\R$.
\end{lemma}

\paragraph{From TFL to Pinsker.} Consider any bounded function $f\colon\domain\to\R$. From~\autoref{lemma:tfl} followed by an application of Hoeffding's Lemma, we can write
\begin{align*}
  \shortexpect_{\p}[f(X)] 
    &\leq \kldiv{\p}{\q} + \log \shortexpect_{\q}[e^{f(Y)}]  \tag{TFL}\\
    &\leq \kldiv{\p}{\q} + \shortexpect_{\q}[f(Y)] + \frac{1}{2}\norminf{f}^2 \tag{Hoeffding's Lemma}
\end{align*}
so, reorganizing and taking the supremum over all $f$ such that $\norminf{f}=\lambda$ (for some $\lambda>0$ to be carefully chosen), we get
\begin{equation}
  \label{eq:from:tfl:to:pinsker}
    \sup_{f: \norminf{f}=\lambda} \mleft( \shortexpect_{\p}[f(X)] - \shortexpect_{\q}[f(Y)] \mright) \leq \kldiv{\p}{\q} + \frac{1}{2}\lambda^2\,.
\end{equation}
Noting then that the LHS is exactly equal to $2\lambda\totalvardist{\p}{\q}$,\footnote{By definition of TV distance as integral probability metric~\cite{Muller97}, or, without using those fancy terms, checking that 
\[
  \totalvardist{\p}{\q} = \sup_{S\subseteq \domain} (\shortexpect_{\p}[\indicSet{S}(X)]-\shortexpect_{\q}[\indicSet{S}(Y)]) = \frac{1}{2}\sup_{f:\norminf{f}\leq 1} (\shortexpect_{\p}[f(X)]-\shortexpect_{\q}[f(Y)])\,.
\]}
from~\eqref{eq:from:tfl:to:pinsker} we are left with the following, true for all $\lambda>0$:
\begin{equation}
  \totalvardist{\p}{\q} \leq \frac{1}{2\lambda}\kldiv{\p}{\q} + \frac{\lambda}{4}\,.
\end{equation}
Optimizing for $\lambda>0$, we choose $\lambda \eqdef \sqrt{2\kldiv{\p}{\q}}$ and obtain
\begin{equation}
  \totalvardist{\p}{\q} \leq \sqrt{\frac{1}{2}\kldiv{\p}{\q}}\,,
\end{equation}
retrieving Pinsker's inequality~\eqref{eq:pinsker}.\smallskip

It is, however, unclear if one can obtain the Bretagnolle--Huber inequality in a similar fashion. At the very least, I do not know how.
\begin{question}
  Can one derive~\eqref{eq:bh} from the TFL?
\end{question}

\paragraph{Update (Aug, 2023):} The answer is yes! Here is a very elegant proof, communicated to me by Hao-Chung Cheng. First, reparameterizing the TFL by setting $g = e^f$, we can write, for any two $\p,\q$,
\[
    \kldiv{\p}{\q} = \sup_{g\colon \Omega\to\R_+} \mleft( \shortexpect_{\p}[\log g(X)] - \log \shortexpect_{\q}[g(Y)] \mright)
\]
and so
\begin{align*}
	e^{-\kldiv{\p}{\q}} &= e^{\inf_{g\colon \Omega\to\R_+} \mleft( \shortexpect_{\p}[-\log g(X)] + \log \shortexpect_{\q}[g(Y)] \mright)}\\
	&= \inf_{g\colon \Omega\to\R_+} e^{\shortexpect_{\p}[-\log g(X)]}\cdot \shortexpect_{\q}[g(Y)] \\
	&\leq \inf_{g\colon \Omega\to\R_+} \shortexpect_{\p}[\tfrac{1}{g(X)}]\cdot \shortexpect_{\q}[g(Y)] \tag{Jensen}
\end{align*}
Since the RHS is an infimum, any choice of $g$ will give an upper bound. So ``all that remains'' to obtain the BH bound~\eqref{eq:bh} is to ``magically'' find a good function $g$ such that $\shortexpect_{\p}[\tfrac{1}{g(X)}]\cdot \shortexpect_{\q}[g(Y)] = 1- \totalvardist{\p}{\q}^2$. Or, factorizing, it would be enough to find $g$ such that $\shortexpect_{\p}[\tfrac{1}{g(X)}] =1 + \totalvardist{\p}{\q}$ and $\shortexpect_{\q}[g(Y)] = 1- \totalvardist{\p}{\q}$. 

Recalling the convenient fact that $\totalvardist{\p}{\q} = 1 - \int_\Omega \min(\p,\q)$, one then chooses 
\[
g(x) = \min\mleft( 1, \frac{\p(x)}{\q(x)}\mright)
\]
(assuming for simplicity here that $\q(x),\p(x)>0$ for all $x$), and\dots{} \emph{that just works}. This gives
\begin{align*}
	e^{-\kldiv{\p}{\q}} &\leq 
	\mleft(\int \max(\p,\q)\mright) \mleft(\int \min(\p,\q)\mright) \\
	&= \mleft(\int \mleft(\p+\q - \min(\p,\q)\mright)\mright) \mleft(\int \min(\p,\q)\mright) \\
	&= \mleft(1 + \totalvardist{\p}{\q}\mright) \mleft(1- \totalvardist{\p}{\q}\mright) \\
	&= 1- \totalvardist{\p}{\q}^2\,
\end{align*}
using that $\p$, $\q$ both integrate to one. Reorganizing the terms, this establishes~\eqref{eq:bh}. $\qed$

\section{Discussion and pointers}
  \label{sec:discussion}
This note is only a succinct, non-exhaustive discussion of possible improvements to Pinsker's inequality, and barely scratches the surface of the many results on this and related questions. We conclude with a few pointers for the interested and fearless reader: Reid and Williamson~\cite{ReidW09} provide a generalization of Pinsker-type inequalities for other $f$-divergences, as well as an  (optimal) integral form of the inequality. Speaking of inequalities between $f$-divergences, Sason and Verd\'u develop in~\cite{SasonV16} techniques to obtain many bounds, among which the Bretagnolle--Huber one. Finally, the recent book of Lattimore and Szepesv{\'a}ri~\cite{LaSze20} covers the BH bound in its chapter on relative entropy (Theorem~14.2), where it also provides some context and discussion.

We could not conclude without mentioning that the excellent lecture notes of Yihong Wu~\cite{Wu20} devote an entire chapter (Section~5) to inequalities between $f$-divergences, including a wonderful theorem due to Harremo\"es and Vajda (Theorem~5.1), which essentially states that to prove \emph{any} such inequality it suffices to prove it for Bernoulli random variables. Those lecture notes also provide, in Section~5.2.2, a handy (albeit short) discussion of Pinsker's inequality, and states the following improvement due to Vajda~\cite{Vajda70}:
\begin{equation}
  \label{eq:vadja}
  \kldiv{\p}{\q} \geq \log\frac{1+\totalvardist{\p}{\q}}{1-\totalvardist{\p}{\q}} - \frac{2\totalvardist{\p}{\q}}{1+\totalvardist{\p}{\q}}
\end{equation}
This is even tighter than the BH bound (\autoref{theo:bh}) we spent so much time covering here, and which only states that $\kldiv{\p}{\q} \geq \log\frac{1}{1-\totalvardist{\p}{\q}^2}$;\footnote{In particular,~\eqref{eq:vadja} does not lose that asymptotic factor $\sqrt{2}$ over Pinsker's for small KL, unlike the BH bound.} however, it is (at least in my eyes) \emph{much} more cumbersome to use.

\bibliographystyle{alpha}
\bibliography{pinskers-and-beyond}
\end{document}